\documentclass[12pt]{article}
\usepackage{amsmath, amsfonts, amssymb, amsthm}
\usepackage{hyperref}
\usepackage{srcltx}
\usepackage{paralist}
\usepackage{graphicx}
\usepackage{color}

\DeclareMathSymbol{\leqslant}{\mathalpha}{AMSa}{"36} 
\DeclareMathSymbol{\geqslant}{\mathalpha}{AMSa}{"3E} 
\DeclareMathSymbol{\eset}{\mathalpha}{AMSb}{"3F}     
\renewcommand{\leq}{\;\leqslant\;}                   
\renewcommand{\geq}{\;\geqslant\;}                   

\newcommand{\comment}[1]{}

\newcommand{\be}{\begin{eqnarray*}}
\newcommand{\ee}{\end{eqnarray*}}
\newcommand{\ben}{\begin{eqnarray}}
\newcommand{\een}{\end{eqnarray}}

\theoremstyle{plain}

\theoremstyle{definition}

\begin{document}

\vglue20pt \centerline{\Large\bf A path formula for the sock sorting problem}
\medskip

\bigskip
\bigskip


\centerline{Simon Korbel and Peter M\"orters}
\bigskip

\begin{center}\it
Institut f\"ur Mathematik\\
Universit\"at zu K\"oln\\
Weyertal 86-90\\
50931 K\"oln\\
Germany\\
\vspace{0.3cm}

\end{center}

\bigskip

{\leftskip=1truecm
\rightskip=1truecm
\baselineskip=15pt
\small

\noindent{\slshape\bfseries Summary.}  Suppose $n$ different pairs of socks are put in a tumble dryer. When the dryer is finished socks are taken out one by one,
if a sock matches one of the socks on the sorting table both are removed, otherwise it is put on the table until its partner emerges from the dryer. 
We note the number of socks on the table after each of the $2n$ socks is taken from the dryer and give an explicit formula for the probability that this sequence
equals a given sequence of length~$2n$.
\bigskip



}

\section{Background and statement of the result}

Suppose that $n$ different pairs of socks are put in a tumble dryer. After operating the dryer the socks are taken out one by one,
if a sock matches one of the socks on the sorting table both are removed, otherwise it is put on the table until its partner emerges from the dryer. 
Let $X_k, k=1,\ldots,2n$ be the number of socks on the table when the $k$th sock is taken from the dryer. How can we describe this process?
\medskip

This problem was known to Daniel Bernoulli~\cite{Bernoulli} who calculated the expectation of $X_k$.
It is also not hard to calculate the variance of $X_k$ and derive a law of large numbers as $n\to\infty$, see \cite{Korbel}.
Some more work allows the derivation of a functional central limit theorem~\cite{LiPritchard}, and to determine the asymptotic
maximum of the process~\cite{Steinsaltz}. But the focus here is on combinatorial, and in particular limit free, results.
\medskip

An observation one can find in the literature \cite{Gilliand, Pantic} is that the paths of our process $(X_1,\ldots,X_{2n})$  are Dyck paths\footnote{A path $(x_1,\ldots, x_{2n})$ is a \emph{Dyck path} if $x_1=1$, $|x_i-x_{i+1}|=1$, $x_i\geq0$ for $1\leq i \leq 2n-1$ and $x_{2n}=0$.} and therefore can be enumerated 
by the Catalan numbers. This observation, however, is not really useful for the sock sorting problem, as different paths have different probabilities, see for example Table~1. 
We therefore go back to a Laplace experiment, by making all socks distinguishable, and use this to derive a formula for the probability of each path.%
\medskip%

To state our formula, we first reduce our path to length $n$ by only noting the number $K_i, i=1,\ldots,n$ of socks on the table when the $i$th pair is matched.
Note that the path $(X_1,\ldots,X_{2n})$ has $n$ upward and $n$ downward steps. $(K_1,\ldots,K_{n})$ represent the heights of the path from which each downward 
step is taken. The full path can be reconstructed from $(K_1,\ldots,K_{n})$, see Lemma~2 below. Our result is the following.
\bigskip%

{\bf Theorem.} If $(k_1,k_2,\ldots,k_n)\in\mathbb N^n$ satisfies 
\begin{equation}\label{1}
k_{i+1}\geq k_i-1 \mbox{ for $1\leq i <n$ and } k_n=1,
\end{equation}
then 
$${\mathbb P}(K_1=k_1,K_2=k_2,\ldots,K_n=k_n) =2^n  \frac{n! \prod^n_{i=1} k_i} {(2n)!}.$$
Otherwise, if $(k_1,k_2,\ldots,k_n)\in\mathbb N^n$ does not satisfy  \eqref{1}, the probability is zero.
\bigskip

For illustration purposes, here is a table of the resulting probabilities of all paths of length ten, characterised by the 
tuples $(k_1, k_2, k_3, k_4, k_5)$ satisfying~\eqref{1}.
\bigskip

\begin{figure}[h]
\centering
\includegraphics{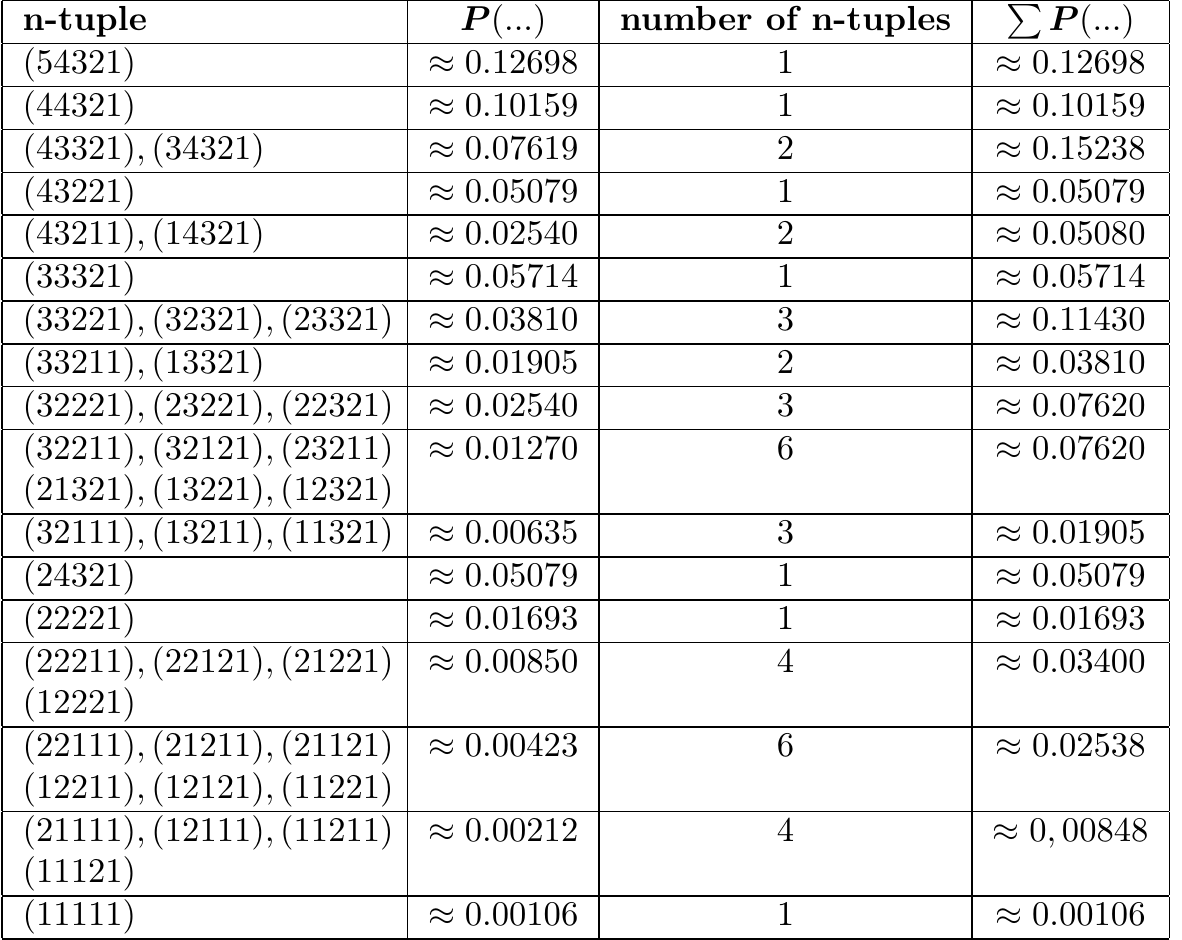}
\caption{Probability of paths of length ten. The probability of a path depends on the tuple
$(k_1, k_2, k_3, k_4, k_5)$ but not on the order of the entries in a tuple.}
\end{figure} 

\section{Derivation of the result}

For a formal definition of $(K_1,\ldots,K_{n})$ we suppose  $x=(x_1,\ldots,x_{2n})$  is a given Dyck path. 
Then we define $L_0(x)=0$ and for $j=1,\ldots,n$ inductively 
$$ L_j(x)= \min\{i > L_{j-1}(x) \colon x_{i+1}= x_i-1\} \mbox{ and } K_j(x)=x_{L_j(x)}. $$
See Figure~2 for an illustration.\\[-2cm]
\begin{figure}[h]
\centering
\includegraphics[width=12cm]{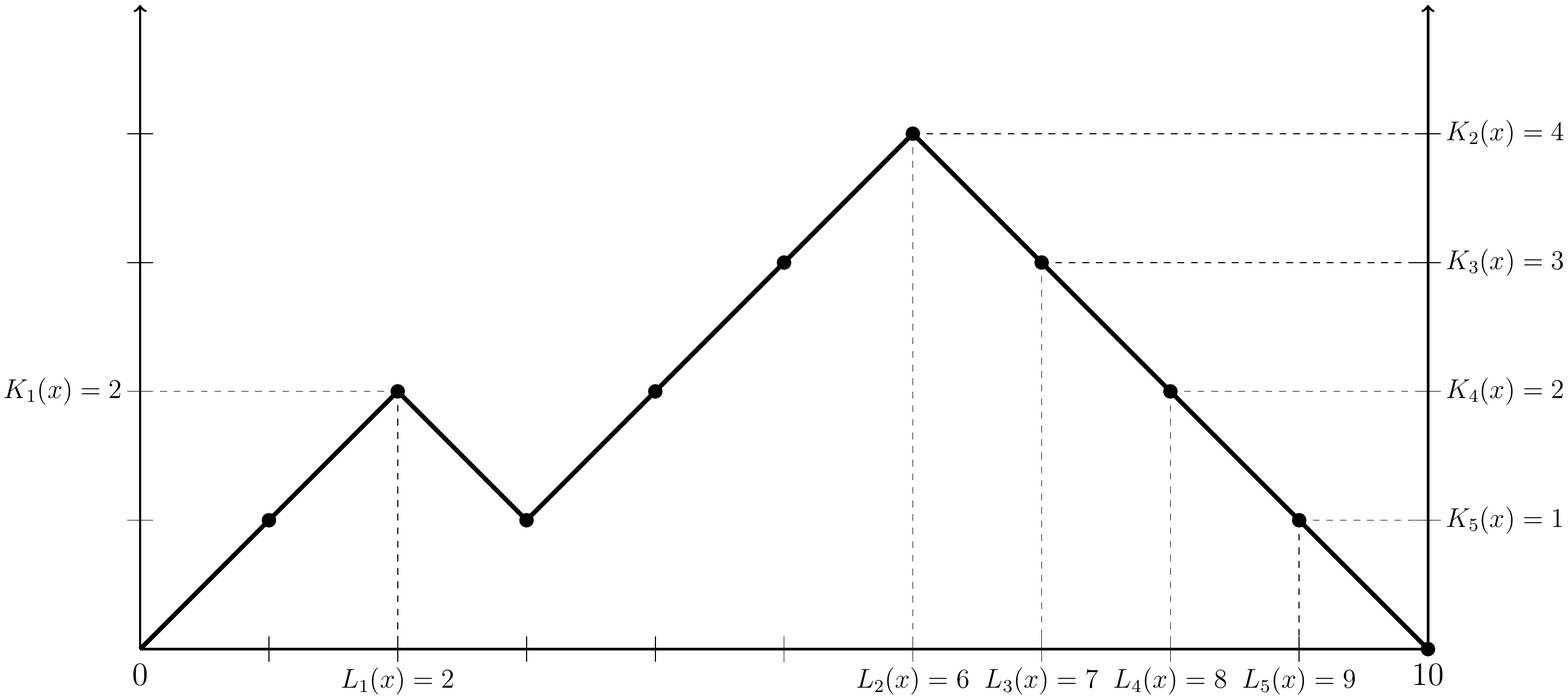}
\\[-1.6cm]
\caption{The path of length ten characterised by the tuple $(2,4,3,2,1)$.}
\end{figure} 

{\bf Lemma 1} If $(x_1,\ldots,x_{2n})$ is a Dyck path and $K_i(x)=k_i$, then 
$(k_1,k_2,\ldots,k_n)$ satisfies~\eqref{1}.%
\bigskip%

{\bf Proof.} Suppose  $(x_1,\ldots,x_{2n})$  is a given Dyck path. Then $L_1(x),\ldots, L_n(x)$ are the ordered elements of the set
$\{i \in\{1,\ldots, 2n-1\} \colon  x_{i+1}= x_i-1\}.$
In particular, as $x_{2n-1}=1$ and $x_{2n}=0$, we have $L_n(x)=2n-1$ and $K_n(x)=1$.\smallskip

Further, by construction, if $L_i(x)<j<L_{i+1}(x)$, then $x_{j+1}=x_j+1$. Hence, 
\begin{align*}
K_{i+1}(x) & =x_{L_{i+1}(x)}= x_{L_{i}(x)} +\sum_{j=L_{i}(x)}^{L_{i+1}(x)-1} (x_{j+1}-x_j)\\ & =K_i(x)-1+(L_{i+1}(x)-L_{i}(x)-1) \geq K_i(x)-1,
\end{align*}
for $i\in\{1,\ldots,n-1\}$. This shows that $(k_1,k_2,\ldots,k_n)$ satisfies~\eqref{1}.
\bigskip

We now show how to reconstruct the Dyck path given $(k_1,k_2,\ldots,k_n)$. Lemma~1 and the construction together establish a bijection 
between the Dyck paths of length $2n$ and the tuples $(k_1,k_2,\ldots,k_n)$ satisfying \eqref{1}. In particular, the cardinality
of the set of such tuples is also given by the Catalan numbers. 
\bigskip

{\bf Lemma 2} If $(k_1,k_2,\ldots,k_n)\in\mathbb N^n$ satisfies 
\eqref{1}, then $x=(x_1,\ldots,x_{2n})$ given by
$$x_i=i  \mbox{ if }i\leq k_1,  \quad 
x_{k_j+2j-1+i}=k_j-1+i \mbox{ if } k_j-1+i\leq k_{j+1}
\mbox{ for $j=1,\ldots,n$.}$$
is the unique Dyck path with $K_i(x)=k_i$.
\pagebreak[3]\bigskip

{\bf Proof.} Suppose  $(k_1,\ldots,k_{n})$ satisfies \eqref{1} and set $k_0=1$. We first check that $x$ as defined in the lemma
is a Dyck path. By definition we have $x_{m+1}-x_m=1$ except when
$m=k_j+2j-1+i$ for $j\in\{0,\ldots, n\}$ and $i\geq 0$ with $k_j-1+i= k_{j+1}$, in which case
$$x_{m+1}=x_{k_j+2j+i}=x_{k_{j+1}+2(j+1)-1}=k_{j+1}-1=k_j-2+i=x_m-1.$$
Hence $x$  has only increments $\pm 1$ and when $x_{m+1}=x_m-1$ we have
$$x_{m+1}=x_m-1=k_{j+1}-1\geq 0.$$
Moreover, as $1\leq k_1$ we have $x_1=1$ and, as $2n=k_n+2n-1$  we have $x_{2n}=k_n-1=0$.
Therefore we have shown that $x$ is a Dyck path. \smallskip

We note from the above that $L_1(x),\ldots, L_n(x)$ are the ordered elements of the set
$$\big\{ m \colon m=k_j+2j-1+i \mbox{ for } i,j \mbox{ with } k_j-1+i= k_{j+1} \big\}.$$
For every $j\in\{0,\ldots, n-1\}$ there is exactly one element in this set which, as
$$k_j+2j-1+i=k_{j+1} +2j < k_{j+1} +2(j+1)-1,$$
is then the  $(j+1)$-smallest element. 
Hence $L_{j+1}(x)=k_j+2j-1+i$  and
$$K_{j+1}(x)=x_{L_{j+1}(x)}=x_{k_j+2j-1+i}=k_j-1+i=k_{j+1}.$$

It remains to show uniqueness. Suppose $x=(x_1,\ldots,x_{2n})$ and
$$\tilde{x}=(x_1,\ldots, x_{m-1} ,\tilde{x}_{m},\ldots, \tilde{x}_{2n})$$ 
are distinct Dyck paths and $m\in\{2,\ldots, 2n-2\}$ is the index of the first step at which they are different.
Without loss of generality we may assume that $x_m=x_{m-1}-1$ and thus
$\tilde{x}_{m}=x_{m-1}+1$. Then there exists $j\in\{1,\ldots,n\}$ with $L_j(x)=m-1$ and $L_j(\tilde{x})>L_j(x)$.
We infer that 
$$K_j(x)=x_{L_j(x)}=x_{m-1} \mbox{ but } K_j(\tilde{x})= x_{m-1}+ (L_j(\tilde{x})-L_j(x))>x_{m-1},$$
showing that $K_j(x)\not=K_j(\tilde{x})$ and thereby implying uniqueness.
\bigskip

%
%
%
%
%
%
%
%

We now look at the probability space $\Omega_n$ consisting of all permutations of $2n$ distinguishable socks. 
In our model each of the $(2n)!$ elements is equally likely. If we write the set of socks as
$$\{1,\ldots, n\} \times \{0,1\},$$
with $0,1$ indicating whether it is the left, resp.\ right, sock,
then $\omega=(\omega_1,\ldots, \omega_{2n})\in\Omega_n$ are the socks in the order drawn from the dryer. Write
$\omega_m=(s_m, p_m)\in\{1,\ldots, n\} \times \{0,1\}$ where
$s_m$ denotes the type of the sock found in the $m$th draw and $p_m$ indicates whether it is the left or right partner. We 
set $s_0=0$ and $x_0=0$ and define
$$x_m=\left\{ \begin{array}{ll} x_{m-1}+1 & \mbox{ if } s_{m}\not\in \{1,\ldots, s_{m-1}\},\\
x_{m-1}-1 & \mbox{ if } s_{m}\in \{1,\ldots, s_{m-1}\},\\
\end{array}\right.$$
and observe that $x(\omega)=(x_1,\ldots, x_{2n})$ is the Dyck path associated with $\omega$.
Our theorem now  follows directly from the following lemma.
\bigskip

{\bf Lemma 3}  If $(k_1,k_2,\ldots,k_n)\in\mathbb N^n$ satisfies 
\eqref{1}, then there are exactly $2^n  n! \prod^n_{i=1} k_i$ different elements $\omega\in\Omega_n$
with $K_i(x(\omega))=k_i$.
\bigskip

{\bf Proof.} Suppose $(k_1,k_2,\ldots,k_n)$ is given and $(x_1,\ldots, x_{2n})$  is the unique Dyck path
with $K_i(x)=k_i$ for all $1\leq i\leq n$. We partition the set of indices $I:=\{1,\ldots,2n\}$ into the sets 
$$M:=\{L_1(x)+1,\ldots, L_n(x)+1\}  \mbox{ and } N:=I \setminus M.$$
At the indices in $N$ the first instance of a sock type is drawn. There are $n!$ ways to allocate
the $n$ types of socks to these $n$ indices and  $2^n$ choices whether the left or right sock is the first.
Once this choice is made we look at the indices in $M$ in order. At each of these indices a pair of socks is completed. 
As there there are $k_i$ socks of different type on the table at time $L_i(x)$, there are exactly $k_i$ choices for 
$\omega_{L_i(x)+1}$. Altogether, we find exactly $2^n  n! \prod^n_{i=1} k_i$ different elements $\omega\in\Omega_n$
with $K_i(x(\omega))=k_i$.
\bigskip

\bigskip

{\bf Acknowledgment:} This note originated from the first author's bachelor thesis.

\end{document}